\documentclass{amsart}
\usepackage{graphicx}
\usepackage[latin2]{inputenc}
\usepackage{amsmath,amssymb,amsbsy,amsthm,amsfonts}

\usepackage{epsfig} 

\begin{document}
\def\note#1{\marginpar{\small #1}}

\def\tens#1{\pmb{\mathsf{#1}}}
\def\vec#1{\boldsymbol{#1}}

\def\norm#1{\left|\!\left| #1 \right|\!\right|}
\def\fnorm#1{|\!| #1 |\!|}
\def\abs#1{\left| #1 \right|}
\def\ti{\text{I}}
\def\tii{\text{I\!I}}
\def\tiii{\text{I\!I\!I}}

\def\diver{\mathop{\mathrm{div}}\nolimits}
\def\grad{\mathop{\mathrm{grad}}\nolimits}
\def\Div{\mathop{\mathrm{Div}}\nolimits}
\def\Grad{\mathop{\mathrm{Grad}}\nolimits}

\def\tr{\mathop{\mathrm{tr}}\nolimits}
\def\cof{\mathop{\mathrm{cof}}\nolimits}
\def\det{\mathop{\mathrm{det}}\nolimits}

\def\lin{\mathop{\mathrm{span}}\nolimits}
\def\pr{\noindent \textbf{Proof: }}
\def\pp#1#2{\frac{\partial #1}{\partial #2}}
\def\dd#1#2{\frac{\d #1}{\d #2}}

\def\T{\mathcal{T}}
\def\R{\mathbb{R}}
\def\bx{\vec{x}}
\def\be{\vec{e}}
\def\bef{\vec{f}}
\def\bec{\vec{c}}
\def\bs{\vec{s}}
\def\ba{\vec{a}}
\def\bn{\vec{n}}
\def\bphi{\vec{\varphi}}
\def\btau{\vec{\tau}}
\def\bc{\vec{c}}
\def\bg{\vec{g}}

\def\bW{\tens{W}}
\def\bT{\tens{T}}
\def\bD{\tens{D}}
\def\bF{\tens{F}}
\def\bB{\tens{B}}
\def\bV{\tens{V}}
\def\bS{\tens{S}}
\def\bI{\tens{I}}
\def\bi{\vec{i}}
\def\bv{\vec{v}}
\def\bfi{\vec{\varphi}}
\def\bk{\vec{k}}
\def\b0{\vec{0}}
\def\bom{\vec{\omega}}
\def\bw{\vec{w}}
\def\p{\pi}
\def\bu{\vec{u}}

\def\ID{\mathcal{I}_{\bD}}
\def\IP{\mathcal{I}_{p}}
\def\Pn{(\mathcal{P})}
\def\Pe{(\mathcal{P}^{\eta})}
\def\Pee{(\mathcal{P}^{\varepsilon, \eta})}

\def\Ln#1{L^{#1}_{\bn}}

\def\Wn#1{W^{1,#1}_{\bn}}

\def\Lnd#1{L^{#1}_{\bn, \diver}}

\def\Wnd#1{W^{1,#1}_{\bn, \diver}}

\def\Wndm#1{W^{-1,#1}_{\bn, \diver}}

\def\Wnm#1{W^{-1,#1}_{\bn}}

\def\Lb#1{L^{#1}(\partial \Omega)}

\def\Lnt#1{L^{#1}_{\bn, \btau}}

\def\Wnt#1{W^{1,#1}_{\bn, \btau}}

\def\Lnd#1{L^{#1}_{\bn, \btau, \diver}}

\def\Wntd#1{W^{1,#1}_{\bn, \btau, \diver}}

\def\Wntdm#1{W^{-1,#1}_{\bn,\btau, \diver}}

\def\Wntm#1{W^{-1,#1}_{\bn, \btau}}

\newtheorem{Theorem}{Theorem}[section]
\newtheorem{Example}{Example}[section]
\newtheorem{Lem}{Lemma}[section]
\newtheorem{Rem}{Remark}[section]
\newtheorem{Def}{Definition}[section]
\newtheorem{Col}{Corollary}[section]
\newtheorem{Proposition}{Proposition}[section]

\newcommand{\Om}{\Omega}
\newcommand{ \vit}{\hbox{\bf u}}
\newcommand{ \Vit}{\hbox{\bf U}}
\newcommand{ \vitm}{\hbox{\bf w}}
\newcommand{ \ra}{\hbox{\bf r}}
\newcommand{ \vittest }{\hbox{\bf v}}
\newcommand{ \wit}{\hbox{\bf w}}
\newcommand{ \fin}{\hfill $\square$}

\newcommand{\ZZ}{\mathbb{Z}}
\newcommand{\CC}{\mathbb{C}}
\newcommand{\NN}{\mathbb{N}}
\newcommand{\V}{\zeta}
\newcommand{\RR}{\mathbb{R}}
\newcommand{\EE}{\varepsilon}
\newcommand{\Lip}{\textnormal{Lip}}
\newcommand{\XX}{X_{t,|\textnormal{D}|}}
\newcommand{\PP}{\mathfrak{p}}
\newcommand{\VV}{\bar{v}_{\nu}}
\newcommand{\QQ}{\mathbb{Q}}
\newcommand{\HH}{\ell}
\newcommand{\MM}{\mathfrak{m}}
\newcommand{\rr}{\mathcal{R}}
\newcommand{\tore}{\mathbb{T}_3}
\newcommand{\Z}{\mathbb{Z}}
\newcommand{\N}{\mathbb{N}}

\newcommand{\F}{\overline{\boldsymbol{\tau} }}

\newcommand{\moy} {\overline {\vit} }
\newcommand{\moys} {\overline {u} }
\newcommand{\mmoy} {\overline {\wit} }
\newcommand{\g} {\nabla }
\newcommand{\G} {\Gamma }
\newcommand{\x} {{\bf x}}
\newcommand{\E} {\varepsilon}
\newcommand{\BEQ} {\begin{equation} }
\newcommand{\EEQ} {\end{equation} }
\makeatletter
\@addtoreset{equation}{section}
\renewcommand{\theequation}{\arabic{section}.\arabic{equation}}

\newcommand{\hs}{{\rm I} \! {\rm H}_s}
\newcommand{\esp} [1] { {\bf I} \! {\bf H}_{#1} }

\newcommand{\vect}[1] { \overrightarrow #1}

\newcommand{\hsd}{{\rm I} \! {\rm H}_{s+2}}

\newcommand{\HS}{{\bf I} \! {\bf H}_s}
\newcommand{\HSD}{{\bf I} \! {\bf H}_{s+2}}

\newcommand{\hh}{{\rm I} \! {\rm H}}
\newcommand{\lp}{{\rm I} \! {\rm L}_p}
\newcommand{\leb}{{\rm I} \! {\rm L}}
\newcommand{\lprime}{{\rm I} \! {\rm L}_{p'}}
\newcommand{\ldeux}{{\rm I} \! {\rm L}_2}
\newcommand{\lun}{{\rm I} \! {\rm L}_1}
\newcommand{\linf}{{\rm I} \! {\rm L}_\infty}
\newcommand{\expk}{e^{ {\rm i} \, {\bf k} \cdot \x}}
\newcommand{\proj}{{\rm I}Ę\! {\rm P}}

\renewcommand{\theenumi}{\Roman{section}.\arabic{enumi}}

\newcounter{taskcounter}[section]

\newcommand{\bib}[1]{\refstepcounter{taskcounter} {\begin{tabular}{ l p{13,5cm}} \hskip -0,2cm [\Roman{section}.\roman{taskcounter}] & {#1}
\end{tabular}}}

\renewcommand{\thetaskcounter}{\Roman{section}.\roman{taskcounter}}

\newcounter{technique}[section]

\renewcommand{\thetechnique}{\roman{section}.\roman{technique}}

\newcommand{\tech}[1]{\refstepcounter{technique} {({\roman{section}.\roman {technique}}) {\rm  #1}}}

\newcommand{\B}{\mathcal{B}}

\newcommand{\diameter}{\operatorname{diameter}}


\title[  Critical and subcritical $\alpha$ models of turbulence]{Mathematical results for some $\displaystyle{\alpha}$ models of turbulence with critical and subcritical regularizations }

\author[H. Ali]{Hani Ali}
\address{IRMAR , UMR CNRS 6625, Universit\'{e} Rennes1, Campus Beaulieu, 35042 Rennes cedex, France}
\email{hani.ali@univ-rennes1.fr}


\keywords{turbulence model, existence, weak solution}
\subjclass[2000]{35Q30,35Q35,76F60}

\begin{abstract}
In this paper, we establish the existence of  a unique  ``regular'' weak   solution  to   turbulent flows governed by a general family of $\alpha$ models with critical regularizations. In particular this family contains  the simplified Bardina  model and the modified Leray-$\alpha$ model. 
When the regularizations are subcritical,
we prove the existence of weak solutions and  we establish an upper bound on the Hausdorff dimension of the time singular set of those weak solutions. 
The result is an  interpolation between the bound proved by Scheffer for the Navier-Stokes equations and the regularity result in the critical case.
\end{abstract}

\maketitle
\maketitle


\section{Introduction}
Let $\mathbb{T}_3 $  be  the  three dimensional torus $\mathbb{T}_3=\left ( \R^3 / {\mathcal T}_3 \right) $ where ${\mathcal T}_{3} = 2 \pi \Z^{3} /L$, $ L > 0$, $0 \le \theta_1, \theta_2 \le {1},$
 and $T\in
(0, \infty)$.
Our goal is to prove, for a given $\bef:(0,T)\times \mathbb{T}_3
\to \mathbb{R}^3$, the existence of  $
(\bv,p):(0,T)\times \mathbb{T}_3 \to \mathbb{R}^3
\times \mathbb{R}$ which solves in a certain sense the following problem $\mathcal{NS}(\alpha) $  
\begin{align}
\diver \bv &=0, \label{BM}\\
\bv_{,t} + \diver (\widetilde{\bv} \otimes \overline{\bv}) -  \nu \Delta \bv  &= -\nabla p + \bef,\label{BLM}\\
\alpha^{2\theta}( -\Delta)^{\theta_1} \widetilde{\bv} +   \widetilde{\bv}
&=\bv, \quad \diver \widetilde{\bv} =0,\label{TKEbis}\\
\alpha^{2\theta}( -\Delta)^{\theta_2} \overline{\bv} +   \overline{\bv}
&=\bv, \quad \diver \overline{\bv} =0.\label{TKE}
\end{align}
considered in $(0,T)\times \mathbb{T}_3$ and completed by appropriate
boundary and initial conditions. Here, $\bv$ is the fluid velocity field, $p$ is  the pressure, $
\bef$ is the  external body forces, $\nu$ stands for the
viscosity.\\
The nonlocal operator $(-\Delta)^{\theta_i}$ , $i=1,2 $ is
defined through the Fourier transform
\begin{equation}
\widehat{(-\Delta)^{\theta_i}{\bv(\bk)}}=|{\bk}|^{2\theta_i}\widehat{\bv}({\bk}).
\end{equation}
Fractionnal order Laplace operator has been used in  another $ \alpha$ models of turbulence in  \cite{OT2007,A01,HLT10}.
Existence and uniqueness of solutions of other modifications of the Navier-Stokes equations have been studied by Ladyzhenskaya \cite{lad67} Lions \cite{Lions59}, M\'{a}lek et al. \cite{mnrr96}.

Our task is to find the critical relation between the  regularizations $\theta_1$ and $\theta_2$ (see Theorem \ref{TH1}) needed to establish global in time existence of a unique weak solution to eqs. \eqref{BM}--\eqref{TKE} and fulfilling the requirements:\\
 $(\bv,p)$ are spatially periodic with period $L$,
\begin{equation}
\label{bc1}
\int_{\tore}\bv(t,\bx) d\bx=0 \quad \textrm{ and } \int_{\tore}p(t,\bx) d\bx=0 \quad \textrm{ for } t \in [0,T),
\end{equation}
and  
\begin{equation}
\bv(0,x)=\bv_0(x) \quad \textrm{ in } \mathbb{T}_3.
\label{ID}
\end{equation}

%
%
%
%
%
Concerning  the regularized velocities $\widetilde{\bv}$, $\overline{\bv}$ we deduce from (\ref{TKEbis}) and (\ref{TKE})  that they verify the same boundary conditions as $\bv$:
 \begin{align}
\widetilde{\bv}(t, \bx + L \vec{e}_{\vec{j}})=\widetilde{\bv}(t,\bx) \quad \textrm{ and  } \int_{\tore}\widetilde{\bv}(t,\bx) d\bx=0 \quad &&\textrm{ on } (0,T)\times \mathbb{T}_3 .\label{bctilde123}\\
\overline{\bv}(t, \bx + L \vec{e}_{\vec{j}})=\overline{\bv}(t,\bx) \quad \textrm{ and  } \int_{\tore}\overline{\bv}(t,\bx) d\bx=0 \quad &&\textrm{ on } (0,T)\times \mathbb{T}_3,\label{bc123}
\end{align}
 We note that the $\alpha$ family considered here is a particular case of the general study  in \cite{HLT10} where the results do not recover the critical case $2\theta_1+\theta_2=\frac{1}{2}$. The  Leray-$\alpha$ model with critical regularization is studied in \cite{A01}.
 We know, thanks to the works \cite{ILT05,CLT06} that for $\theta_1=0,\  \theta_2=1$ or  $\theta_1=1, \ \theta_2=1, $  that their exist a unique weak solution to the model \eqref{BM}--\eqref{TKE}.\\
  When  $\theta_1=1, \ \theta_2=1, $ we get the simplified Bardina model \cite{CLT06}.  The simplified Bardina model first arose in the context of turbulence models for the Navier-Stokes equations in \cite{LL03}. Based on this work, we will study in a forthcoming paper the model studied in  \cite{LL03,LL06b} and  other related model \cite{dunca06,roger2011} in the special case where the filtering is given by 
$ \alpha^{2\theta}( -\Delta)^{\theta} \overline{\phi} +   \overline{\phi}
=\phi.$ \\
 When the relation between the  regularizations $\theta_1$ and $\theta_2$ is subcritical 
 we will prove that $\frac{1-2\theta_2-4\theta_1}{2}$-dimensional Hausdorff measure of  the time singular set $\mathcal{S}_{\theta_1, \theta_2}(\overline{\vec{v}})$ of any 
weak solution $\overline{\vec{v}}$ of  \eqref{BM}--\eqref{TKE} is zero (see Theorem \ref{main}).\\
The Hausdorff dimension of the time singular set to weak solutions of another modification of the Navier-Stokes equations was studied in \cite{BumaPr11,AA11}.\\
As a conclusion our study gives the critical regularizations  to various $\alpha$ models, namley  the modified Leray-$\alpha$ \cite{ILT05} and the simplified Bardina model \cite{CLT06}.
These critical regularisations and the Hausdorff measure of the time singular set in the subcritical case are listed in table \ref{Table1}.
 \begin{table}[h]
\begin{tabular}{|c|c|c|c|c|r|}
  \hline
   &  \small{simplified Bardina}& \small{Leray}-$\alpha$  & \small{modified Leray-}$\alpha$  \\
  \hline
   \hline
   
$  \displaystyle \theta_1$ &$ \displaystyle \frac{1}{6} $ &$ \displaystyle \frac{1}{4}$ &$\displaystyle 0$  \\

  &   &  &  \\
 $ \displaystyle \theta_2$     & $\displaystyle \frac{1}{6}$     & $\displaystyle 0$                      &$\displaystyle \frac{1}{2}$  \\
  &   &  &  \\
  $ \displaystyle \mathcal{H}(\mathcal{S})$ &    $\displaystyle \frac{1-6\theta_1}{2} $     &$\displaystyle \frac{1-4\theta_1}{2} $   
  & $\displaystyle \frac{1-2\theta_2}{2}$     \\
 \hline
    \end{tabular}
    \caption{\small{Comparison of various critical regularizations and Hausdorff measure  for the simplified Bardina, Leray-$\alpha$ and modified Leray-$\alpha$}}
    \label{Table1}
   \end{table}\\
   Observe that the results reported here are also valid in the whole space $\R^3$ by
employing the relevant analogue tools for treating the Navier-Stokes in the whole
space \cite{shonbek,alihmidi}.

This  paper is organized as follows. Section 2 consists of notation and conventions used throughout. In section 3 we prove the global existence and uniqueness of the solution to the 
 model \eqref{BM}--\eqref{TKE} with critical regularization.  Section 4 treats the
question of the subcritical regularizations where we give  an upper bound on the Hausdorff dimension of the time singular set of  weak solutions to   the model \eqref{BM}--\eqref{TKE}.   
The result is an  interpolation between the bound proved by Scheffer for the Navier-Stokes equations and the regularity result in the critical case.
   
   \section{Notations}
\ Before formulating the main results of this paper, we fix notation of function spaces that we shall  employ.\\ 
We denote by $L^p(\tore)$ and $W^{r,p}(\tore)$, $r \ge -1, \ 1 \le p \le \infty,$ the usual Lebesgue and Sobolev spaces over $\tore$, and the Bochner spaces $C(0,T;X), L^p(0,T;X)$ are defined in the standard way. 

The Sobolev spaces $ \vec{H}^{s}=H^s(\tore)^3$,  of mean-free functions are classically characterized in terms of
the Fourier series
$$\vec{H}^{s} = \left\{ \bv(\bx)=\sum_{{\vec{k} \in {\mathcal T}_{3}}}^{}\bc_{\vec{k}} e^{i\vec{k}\cdot \vec{x}}, \left(\bc_{\vec{k}} \right)^{*}=\bc_{-\vec{k}}, \bc_0 = 0,   \|\bv \|_{s,2}^2= \sum_{{\vec{k} \in {\mathcal T}_{3}}}^{}| \vec{k} |^{2s} |\bc_{\vec{k}}|^2<\infty  \right\},$$
where  $ \left(\bc_{\vec{k}}^N \right)^{*}$ denote the complex conjugate $\bc_{\vec{k}}^N.  $
In addition we introduce
\begin{align*}
\vec{H}^{s}_{\diver}&=\left\{ \bv \in \vec{H}^s; \; \diver \bv
=0 \textrm{ in }  \tore \right\},\\
\vec{H}^{-s}&=\left(\vec{H}^{s}_{} \right)^{'},\quad
\vec{L}^2_{}=\vec{H}^{0},\quad 
\vec{L}^2_{\diver}=\vec{H}^{0}_{\diver}.
\end{align*}
 Let us mention that by using Poincar\'e inequality we have 
  \begin{align}
 \|\vec{v} \|_{s,2} \approx \| \widetilde{\vec{v}}\|_{s+2\theta_1,2} \approx \| \overline{\vec{v}} \|_{s+2\theta_2,2}. 
 \end{align}
Throughout we will use $C$ to denote an arbitrary constant which may change line to
line.
\section{Existence and uniqueness in the critical case: $2\theta_1 +   \theta_2 = \frac{1}{2} $}
 The aim in this section is to find the critical relation between $\theta_1$ and $\theta_2$ that ensures  the existence and the uniqueness of the weak solution to the model \eqref{BM}--\eqref{TKE}. 

\begin{Theorem}
Assume that $2\theta_1 + \theta_2 = \frac{1}{2}$, $0 \le \theta_1 < \frac{1}{4}$ and $0 < \theta_2 \le \frac{1}{2}$. Let  $\bef \in L^{2}
(0,T;\vec{H}^{-1}_{})$ be a divergence free function and  $\bv_0 \in L^2_{\diver}$. Then there exist $(\bv, p)$ a unique ``regular'' weak  solution to \eqref{BM}--\eqref{TKE}  such that
\begin{align}
\bv &\in \mathcal{C}_{}(0,T;\vec{L}^2_{\diver}) \cap
L^2(0,T;\vec{H}^{1}_{\diver}),\label{bv12}\\
\bv_{,t}&\in   L^{2}(0,T;\vec{H}^{-1}_{}),
\label{bvt}\\
p&\in L^{2}(0,T;L^{2}(\tore))
\label{psp}.
\end{align}
fulfill
\begin{equation}
\begin{split}
\int_0^T \langle \bv_{,t}, \bw \rangle  -  (\widetilde{\bv} \otimes \overline{\bv}, \nabla
\bw) +  \nu(
\nabla \bv, \nabla \bw  )\; dt
-(p, \diver \bw) \\
=\int_0^T  \langle \bef, \bw \rangle \; dt
\qquad \textrm{ for all } \bw\in L^{2}(0,T; \vec{H}^{1}_{}),
\end{split}\label{weak1}
\end{equation}
Moreover,
\begin{equation}
\bv(0)= \bv_0.
\label{intiale}
\end{equation}
 \label{TH1}
\end{Theorem}

\begin{Rem}
We use the name ``regular"  for the weak  solution since the weak solution is unique and the velocity part of the solution $\bv$ is a possible test function in the weak formulation (\ref{weak1}), that in particular implies that  $\bv \in  \mathcal{C}_{}(0,T;\vec{L}^2_{\diver}).$
\end{Rem}
\begin{Rem}
 Once existence and  uniqueness in the large of a weak solution to the model \eqref{BM}--\eqref{TKE} with critical regularization is known. Further theoretical
properties of the model  can then be developed. These are
currently under study by the author and will be presented in a subsequent report.
\end{Rem}
\begin{Rem}
In a further paper, the author will prove that the solution $(\overline{\bv},p)$ of the model \eqref{BM}--\eqref{TKE}   converges in
some sense to a solution of the Navier-Stokes equations  when $\alpha$  goes to zero. 
\end{Rem}

\textbf{Proof of Theorem \ref{TH1}}\\
The proof of Theorem  \ref{TH1} follows the classical scheme. We start by constructing  approximated solutions $(\bv^N,p^N)$ via Galerkin method. Then we   seek for a priori estimates that are uniform with respect to $N$. Next, we  passe to the limit in the equations after
having used compactness properties. Finaly we show that the solution we constructed is unique thanks to  Gronwall's
lemma. We also note that in our argument we keep 
the pressure in the weak formulation of the problem and we do not simply neglect it by
projecting the equations over divergence-free vector fields.\\

\textbf{Step 1}(Galerkin approximation).
Consider the sequence $\left\{ e^{i\vec{k}\cdot \vec{x}} \right\}_{|\vec{k}|=1}^{\infty}$ consisting of $L^2$-orthonormal  and $W^{1,2}$-orthogonal eigenvectors of the following problem:
\begin{align}
-\Delta e^{i\vec{k}\cdot \vec{x}} = |\vec{k}|^2 e^{i\vec{k}\cdot \vec{x}},\hbox{ in } \tore, \quad  \hbox{ for all } {\vec{k}} \in {\mathcal T}_{3}\setminus\{0\} .
\end{align}
We note that this sequence forms a hilbertian basis of $\vec{L}^2$.\\
We set 
\begin{equation}
\bv^N(t,\bx)=\sum_{{|\vec{k}|=1}}^{N}\bc_{\vec{k}}^N (t) e^{i\vec{k}\cdot \vec{x}}, 
\end{equation}
 such that $\bk \cdot  \bc_{\vec{k}}^N = 0  $ for all $\bk \in {\mathcal T}_{3}\setminus\{0\}$ and $\left(\bc_{\vec{k}}^N \right)^{*}=\bc_{-\vec{k}}^N  $. 
 Thus due of (\ref{TKEbis}) and (\ref{TKE}) we have 
\begin{equation}
 \widetilde{\bv}^N(t,\bx)=\sum_{{|\vec{k}|=1}}^{N}\widetilde{\bc}_{\vec{k}}^N (t) e^{i\vec{k}\cdot \vec{x}} \textrm{ and }  \overline{\bv}^N(t,\bx)=\sum_{{|\vec{k}|=1}}^{N}\overline{\bc}_{\vec{k}}^N (t) e^{i\vec{k}\cdot \vec{x}}, 
\end{equation}
 where 
\begin{equation}
 \widetilde{\bc}_{\vec{k}}^N = \frac{{\bc}_{\vec{k}}^N }{1+\alpha^{2\theta_1}|\vec{k}|^{2\theta_1}} \textrm{ and }  \overline{\bc}_{\vec{k}}^N = \frac{{\bc}_{\vec{k}}^N }{1+\alpha^{2\theta_2}|\vec{k}|^{2\theta_2}}, 
\end{equation}
for all $\bk \in {\mathcal T}_{3}\setminus\{0\}$.\\
We look for $(\bv^N(t,\bx), p^N(t,\bx)) $ that are determined through the system of equations 
\begin{equation}
\begin{split}
 \left( \bv_{,t}^N, e^{i\vec{k}\cdot \vec{x}} \right)  -  (\widetilde{\bv}^N \otimes \overline{\bv}^N, \nabla
e^{i\vec{k}\cdot \vec{x}}) +  \nu(
\nabla \bv^N, \nabla e^{i\vec{k}\cdot \vec{x}}  )\; 
-(p^N, \diver e^{i\vec{k}\cdot \vec{x}})\\
 \qquad = \langle \bef, e^{i\vec{k}\cdot \vec{x}} \rangle \; ,
\qquad {|\vec{k}|=1},2,...,N,
\end{split}\label{weak1galerkine}
\end{equation}
and 
\begin{equation}
\begin{split}
 p^N = - \sum_{i,j}\partial_i\partial_j \Delta^{-1}(\Pi^N(\widetilde{v}_i^N\overline{v}_j^N))= - \sum_{i,j} R_{ij}(\Pi^N(\widetilde{v}_i^N\overline{v}_j^N)).
\end{split}\label{pressuregalerkine}
\end{equation}
Where the projector $ \displaystyle \Pi^N $ assign to any Fourier series $\displaystyle \sum_{\bk \in {\mathcal T}_{3}\setminus\{0\}} \vec{g}_{\bk} e^{i\bk \cdot\bx} $ its N-dimensional part, i.e. 
 $\displaystyle \sum_{\bk \in {\mathcal T}_{3}\setminus\{0\}, |\bk| \le N} \vec{g}_{\bk} e^{i\bk \cdot\bx}, $
 and $R_{ij}$ is the Riez operator defined through the Fourier transform by
 \begin{equation}
\begin{split}
  \widehat{R_{ij}(u)}=\frac{k_i k_j}{|\bk|^2}\widehat{u({\bk})}, \quad \hbox{ for all } \bk \in {\mathcal T}_{3}\setminus\{0\}.
  \end{split}\label{pressurepseudo2}
\end{equation}

Moreover we require that $\bv^N $ satisfies the following initial condition
\begin{equation}
 \label{initial Galerkine}
\bv^N(0,.)= \bv^N_0= \sum_{|\vec{k}|=1}^{N}\bc_0^N  e^{i\vec{k}\cdot \vec{x}},
\end{equation}
and 
\begin{equation}
 \label{initial2 Galerkine}
\bv^N_0 \rightarrow \bv_0  \quad \textrm{ strongly  in } L^{2}(0,T; \vec{L}^2) \quad \textrm{ when } N \rightarrow \infty .
\end{equation}
Where the initial condition $ \overline{\bv}^N_0$ is deduced from $\bv^N_0$ through the relation (\ref{TKE}).

The classical Caratheodory theory \cite{Wa70}  then implies the short-time existence of solutions 
to (\ref{weak1galerkine})-(\ref{pressuregalerkine}).  Next we derive  estimate on $\bc^N$ that is uniform w.r.t. $N$.
These estimates then imply that the  solution of  (\ref{weak1galerkine})-(\ref{pressuregalerkine}) constructed on a short time interval $[0, T^N[ $ exists for all $t \in [0, T]$.\\

\textbf{Step 2} (Uniform estimates 1)
Multilplying the $|\vec{k}|$th equation in (\ref{weak1galerkine}) with $\overline{\bc}^N_{\vec{k}}(t)$, summing over ${|\vec{k}|=1},2,...,N$, integrating over time from $0$ to $t$ and using the following identities
 
\begin{equation}
 \label{divergencfreebar1}
\left({\bv}^N_{,t},  \overline{\bv}^N \right)=\left(\overline{\bv}^N_{,t}+ \alpha^{2\theta}(-\Delta)^\theta\overline{\bv}^N_{,t} ,  \overline{\bv}^N \right)=
\frac{1}{2}\frac{d}{dt}\|\overline{\bv}^N \|_{2}^2 + \frac{1}{2}\frac{d}{dt}\|\overline{\bv}^N \|_{\theta,2}^2,
\end{equation}
\begin{equation}
 \label{divergencfreebar2}
\left(-\Delta{\bv}^N,  \overline{\bv}^N \right)=\left(-\Delta \overline{\bv}^N -\alpha^{2\theta}\Delta(-\Delta)^\theta \overline{\bv}^N ,  \overline{\bv}^N \right)=
\|\overline{\bv}^N \|_{1,2}^2 + \|\overline{\bv}^N \|_{1+\theta,2}^2,
\end{equation}
and
\begin{equation}
 \label{divergencfreebar}
\left(\widetilde{\bv}^N \otimes \overline{\bv}^N, \nabla \overline{\bv}^N \right)=\left(  \widetilde{\bv}^N, \nabla \frac{|\overline{\bv}^N|^2}{2}\right)
=-\left(  \diver \widetilde{\bv}^N,  \frac{|\overline{\bv}^N|^2}{2}\right)=0
\end{equation}
leads to the a priori estimates
\begin{equation}
\begin{array}{llll}
 \label{apriori1}
\displaystyle \frac{1}{2}\left(\|\overline{\bv}^N \|_{2}^2 + \|\overline{\bv}^N \|_{\theta,2}^2\right)+ \displaystyle \nu\int_{0}^{t}\left(\|\overline{\bv}^N \|_{1,2}^2  + \|\overline{\bv}^N \|_{1+\theta_2,2}^2 \right) \ ds\\
 \quad \quad \quad = \displaystyle \int_{0}^{t} \langle \bef, \overline{\bv}^N \rangle \ ds 
 +\displaystyle \frac{1}{2}\left(\|\overline{\bv}_0 \|_{2}^2 + \|\overline{\bv}_0 \|_{\theta_2,2}^2\right).
\end{array}
\end{equation}
 Using the duality norm comined with Young inequality we conclude from eqs. (\ref{apriori1}) that 
 \begin{equation}
 \label{apriori12}
\sup_{t \in [0,T^N[}\|\overline{\bv}^N \|_{2}^2 + \sup_{t \in [0,T^N[}\|\overline{\bv}^N \|_{\theta_2,2}^2+ \nu\int_{0}^{t}\left(\|\overline{\bv}^N \|_{1,2}^2  + \|\overline{\bv}^N \|_{1+\theta_2,2}^2 \right) \ ds  \le C
\end{equation}
that immediately implies that the existence time is independent of $N$ and it is possible to take $T=T^N$.\\ 
We deduce from \ref{apriori12} that 
\begin{equation}
\label{vbar1}
 \overline{\bv}^N \in L^{\infty}(0,T ; \vec{H}^{\theta_2}_{\diver}) \cap L^{2}(0,T ; \vec{H}^{1+\theta_2}(\tore)^3), 
 \end{equation}
 thus from the relation (\ref{TKE}) combined with  the Poincar\'e inequality we conclude that 
  \begin{equation}
\label{v1}
 {\bv}^N \in L^{\infty}(0,T ; \vec{H}^{-\theta_2}_{\diver}) \cap L^{2}(0,T ; \vec{H}^{1-\theta_2}).  
 \end{equation}
 From (\ref{TKEbis}) it follows that 
  \begin{equation}
\label{vtilde1}
 \widetilde{\bv}^N \in L^{\infty}(0,T ; \vec{H}^{2\theta_1 -\theta_2}_{\diver}) \cap L^{2}(0,T ; \vec{H}^{1+2\theta_1-\theta_2}).  
 \end{equation}
 
\textbf{Step 3} (Uniform  estimates 2) Let us come back to the relation  (\ref{weak1galerkine}), multilplying the $|\vec{k}|$th equation in (\ref{weak1galerkine}) with ${\bc}^N_{\vec{k}}(t)$, summing over ${|\vec{k}|=1},2,...,N$,
we conclude that 
\begin{equation}
 \label{apriori2}
\displaystyle  \frac{1}{2}\frac{d}{dt}\|{\bv}^N \|_{2}^2 + \nu \|{\bv}^N \|_{1,2}^2    \le  \left| \left(\widetilde{\bv}^N \otimes \overline{\bv}^N, \nabla {\bv}^N \right) \right|+  \langle \bef, {\bv}^N \rangle  \quad := I_1 + I_2.
\end{equation}
For $I_1$ we have for $ \frac{1}{2}-2\theta_1+\theta_2 \le 2\theta_2$ i.e. $ 2\theta_1+\theta_2 \ge \frac{1}{2}$ that 
 \begin{equation}
 \begin{array}{llll}
 \label{apriori22}
I_1 &\le \|\widetilde{\bv}^N \otimes \overline{\bv}^N \|_{2} \| \nabla {\bv}^N \|_{2}\\
 &\le \displaystyle \frac{C}{\nu}\|\widetilde{\bv}^N \otimes \overline{\bv}^N \|_{2}^2+ \frac{\nu}{4} \| \nabla {\bv}^N \|_{2}^2\\
 &\le\displaystyle \frac{C}{\nu}\|\widetilde{\bv}^N \|_{1+2\theta_1-\theta_2,2}^2 \| \overline{\bv}^N \|_{\frac{1}{2}-2\theta_1+\theta_2,2}^2 + \frac{\nu}{4} \| \nabla {\bv}^N \|_{2}^2\\
 &\le \displaystyle \frac{C}{\nu}\|\widetilde{\bv}^N \|_{1+2\theta_1-\theta_2,2}^2 \|\overline{\bv}^N \|_{2\theta_2,2}^2 + \frac{\nu}{4} \| \nabla {\bv}^N \|_{2}^2.
 \end{array}
 \end{equation}
 Now we use the following inequality (see in \cite{A01}.)
 \begin{equation}
 \label{regularity2teta}
\|\overline{\bv}^N \|_{2\theta_2,2}^2 \le \frac{1}{\alpha^{2\theta_2}}\|{\bv}^N \|_{2}^2.
\end{equation}
We conclude  that 
  \begin{equation}
 \begin{array}{llll}
 \label{apriori222}
 I_1 \le \displaystyle \frac{C}{\nu}  \frac{1}{\alpha^{2\theta}}\|{\bv}^N \|_{2}^2   \|\widetilde{\bv}^N \|_{1+2\theta_1-\theta_2,2}^2 +  \frac{\nu}{4} \| \nabla {\bv}^N \|_{2}^2.
 \end{array}
 \end{equation}
 To estimate $I_2$ we use the duality norm and Young inequality in order to obtain
  \begin{equation}
\label{apriori2222}
I_2 \le \| \bef\|_{-1,2} \|{\bv}^N \|_{1,2} \le \frac{C}{\nu} \| \bef\|_{-1,2}^2 + \frac{\nu}{4} \|  {\bv}^N \|_{1,2}^2.
\end{equation}
 Thus  (\ref{apriori222}) and (\ref{apriori2222})  lead to the conclusion that 
    \begin{equation}
 \begin{array}{llll}
 \label{apriori222222}
 \displaystyle  \frac{1}{2} \frac{d}{dt}\|{\bv}^N \|_{2}^2 + \nu \|{\bv}^N \|_{1,2}^2  \le \displaystyle \frac{C}{\nu}  \frac{1}{\alpha^{2\theta}}\|{\bv}^N \|_{2}^2   \|\widetilde{\bv}^N \|_{1+2\theta_1-\theta_2,2}^2 +  \frac{C}{\nu} \| \bef\|_{-1,2}^2.
 \end{array}
 \end{equation}
 Integrating (\ref{apriori222222}) over time from $0$ to $T$ and using Gronwall's Lemma and (\ref{vtilde1}) lead to the following estimate
 \begin{equation}
 \label{apriori2220}
\sup_{t \in [0,T]}\|{\bv}^N \|_{2}^2 + \nu\int_{0}^{T}\|{\bv}^N \|_{1,2}^2  \ dt  \le C.
\end{equation}
 We deduce from (\ref{apriori2220}) that 
\begin{equation}
\label{v11}
 {\bv}^N \in L^{\infty}(0,T ; \vec{L}^{2}_{\diver}) \cap L^{2}(0,T ; \vec{H}^{1}_{\diver}), 
 \end{equation}
  thus from the relation (\ref{TKE})  we conclude that 
  \begin{equation}
\label{vbar11}
 \overline{\bv}^N \in L^{\infty}(0,T ; \vec{H}^{2\theta_2}_{\diver}) \cap L^{2}(0,T ; \vec{H}^{1+2\theta_2}_{\diver}),
 \end{equation}
 and from (\ref{TKEbis}) we obtain
   \begin{equation}
\label{vtilde11}
 \widetilde{\bv}^N \in L^{\infty}(0,T ; \vec{H}^{2\theta_1}_{\diver}) \cap L^{2}(0,T ; \vec{H}^{1+2\theta_1}_{\diver}).
 \end{equation}
 We observe from (\ref{vbar11}) and (\ref{vtilde11}) that 
{ for all } $2\theta_1 + 2\theta_2 \ge \frac{1}{2}$, {in particular for } $2\theta_1  + \theta_2 \ge \frac{1}{2}$, { we have }
 \begin{equation}
\label{vbarvbar}
\widetilde{\bv}^N \otimes \overline{\bv}^N \in  L^{2}(0,T ;L^{2}(\tore)^{3 \times 3} ). 
\end{equation}
Consequently from the Calderon-Zygmund theory eqs (\ref{pressuregalerkine}) implies that 
 \begin{equation}
\label{vbarvbarpressure}
\int_{0}^{T}\|p^N\|_{2}^2 dt < K. 
\end{equation}
From eqs. (\ref{weak1galerkine}), (\ref{vbar11}) and (\ref{vtilde11}) we also obtain that 
 \begin{equation}
\label{vtemps}
\int_{0}^{T}  \|\bv^N_{,t}\|_{-1,2}^2  dt < K. 
\end{equation}
and thus from the relations (\ref{TKEbis}) and (\ref{TKE}) we deduce 
 \begin{equation}
\label{vtildebartemps}
\int_{0}^{T}  \|\widetilde{\bv}^N_{,t}\|_{-1,2}^2  dt < K, \quad \textrm{ and } \int_{0}^{T}  \|\overline{\bv}^N_{,t}\|_{-1,2}^2 \le K.
\end{equation}

\textbf{Step 4} (Limit $N \rightarrow \infty$) It follows from the estimates (\ref{v11})-(\ref{vtildebartemps}) and the Aubin-Lions compactness lemma
(see \cite{sim87} for example) that there are a  not relabeled  subsequence of $(\bv^N,\widetilde{\bv}^N,\overline{\bv}^N, p^N)$  and a quadruplet $(\bv,\widetilde{\bv},\overline{\bv}, p)$ such that

\begin{align}
\bv^N &\rightharpoonup^* \bv &&\textrm{weakly$^*$ in } L^{\infty}
(0,T;\vec{L}^2_{}), \label{c122}\\
\widetilde{\bv}^N &\rightharpoonup^* \widetilde{\bv} &&\textrm{weakly$^*$ in } L^{\infty}
(0,T;\vec{H}^{2\theta_1}_{}), \label{c122primetilde}\\
\overline{\bv}^N &\rightharpoonup^* \overline{\bv} &&\textrm{weakly$^*$ in } L^{\infty}
(0,T;\vec{H}^{2\theta_2}_{}), \label{c122prime}\\
\bv^N &\rightharpoonup \bv &&\textrm{weakly in }
L^2(0,T;\vec{H}^{1}_{}), \label{c22}\\
\widetilde{\bv}^N &\rightharpoonup \widetilde{\bv} &&\textrm{weakly in }
L^2(0,T;\vec{H}^{1+2\theta_1}_{}), \label{c22primetilde}\\
\overline{\bv}^N &\rightharpoonup \overline{\bv} &&\textrm{weakly in }
L^2(0,T;\vec{H}^{1+2\theta_2}_{}), \label{c22prime}\\
\bv^N_{,t}&\rightharpoonup \bv_{,t} &&\textrm{weakly in } L^{2}
(0,T;\vec{H}^{-1}_{}),
\label{c322}\\
\widetilde{\bv}^N_{,t}&\rightharpoonup \widetilde{\bv}_{,t} &&\textrm{weakly in } L^{2}
(0,T;\vec{H}^{-1}_{}),
\label{c322bistilde}\\
\overline{\bv}^N_{,t}&\rightharpoonup \overline{\bv}_{,t} &&\textrm{weakly in } L^{2}
(0,T;\vec{H}^{-1}_{}),
\label{c322bis}\\
p^N&\rightharpoonup p &&\textrm{weakly in } L^{2}(0,T;L^{2}
(\tore)), \label{c32}\\
\bv^N &\rightarrow \bv &&\textrm{strongly in  }
L^2(0,T;\vec{L}^2),\label{c83ici}\\
\widetilde{\bv}^N &\rightarrow \widetilde{\bv} &&\textrm{strongly in  }
L^2(0,T;\vec{L}^2),\label{c82in}\\
\overline{\bv}^N &\rightarrow \overline{\bv} &&\textrm{strongly in  }
L^2(0,T;\vec{L}^2).\label{c82int}
\end{align}

By a standard interpolation argument we have  

\begin{align}
\bv^N & \in  L^{\frac{10}{3}}(0,T;L^{\frac{10}
{3}}(\tore)^3), \label{intrpolv}\\
\widetilde{\bv}^N & \in  L^{\frac{10}{3-4\theta_1}}(0,T;L^{\frac{10}
{3-4\theta_1}}(\tore)^3), \label{interpvtilde}\\
\overline{\bv}^N & \in  L^{\frac{10}{3-4\theta_2}}(0,T;L^{\frac{10}
{3-4\theta_2}}(\tore)^3). \label{interpvbar}
\end{align}
Thus  from (\ref{intrpolv})-(\ref{interpvbar}) and (\ref{c83ici})-(\ref{c82int})
we obatin
\begin{align}
\bv^N &\rightarrow \bv &&\textrm{strongly in  }
L^{q_1}(0,T;L^{q_1}(\tore)^3) \textrm{ for all } q_1<\frac{10}{3},\label{convint1}\\
\widetilde{\bv}^N &\rightarrow \widetilde{\bv} &&\textrm{strongly in  }
L^{q_2}(0,T;L^{q_2}(\tore)^3) \textrm{ for all } q_2<\frac{10}{3-4\theta_1},\label{convint2}\\
\overline{\bv}^N &\rightarrow \overline{\bv} &&\textrm{strongly in  }
L^{q_3}(0,T;L^{q_3}(\tore)^3) \textrm{ for all } q_3<\frac{10}{3-4\theta_2},\label{convint3}
\end{align}
Since $ q_2<\frac{10}{3-4\theta_1}$, $q_3<\frac{10}{3-4\theta_1} $ and $2\theta_1 + \theta_2 = \frac{1}{2} $  the application of H$\ddot{o}$lder's inequality implies that
\begin{align}
 \widetilde{\bv} \otimes \overline{\bv} \in 
L^{q}(0,T;L^{q}(\tore)^{3\times 3}) \textrm{ where } q \ge {2},\label{vtimesv}
\end{align}

The above established convergences are clearly sufficient for taking the limit in (\ref{weak1galerkine})  and for concluding that   $ (\bv,p)$  satisfy (\ref{weak1}). 
Moreover, 
from (\ref{c22}) and (\ref{c322})
 one we can deduce by a classical argument of J.L. Lions \cite{JL69}  that 
 \begin{equation}
 \bv \in  \mathcal{C}(0,T;\vec{L}^2_{}).
\end{equation}
Furthermore, from  the  strong continuty of $\bv$ with respect to the time with value in $L^2_{}$  we deduce   that $\bv(0)=\bv_0$.\\
Let us mention also that $\overline{\bv} $ is a possible  test in the weak formlation (\ref{weak1}). Thus $ \overline{\bv}$ verifies for all $t \in [0,T]$  the follwing equality   
\begin{equation}
\begin{array}{lll}
 \label{apriori12leray}
\displaystyle \left(\|\overline{\bv} (t)\|_{2}^2 + \|\overline{\bv}(t)\|_{\theta_2,2}^2\right)+\displaystyle 2\nu\int_{0}^{t}\left(\|\overline{\bv} \|_{1,2}^2  + \|\overline{\bv} \|_{1+\theta_2,2}^2 \right)ds \\
\quad \quad \quad = \displaystyle 2\int_{0}^{t}\langle \bef, \overline{\bv} \rangle ds + \left(\|\overline{\bv}_0 \|_{2}^2 + \|\overline{\bv}_0\|_{\theta_2,2}^2\right).
\end{array}
\end{equation}

 \textbf{Step 5} (Uniqueness)
Since the pressure part of the solution is uniquely determined by the velocity part it remain to show the uniqueness to the velocity.

Next, we will show the continuous dependence of the  solutions on the initial data and in particular the uniqueness.\\
Let $({\bv_1,p_1})$ and $({\bv_2,p_2})$   any two solutions of (\ref{BM})-(\ref{TKE}) on the interval $[0,T]$, with initial values $\bv_1(0)$ and $\bv_2(0)$. Let us denote by  $\vec{w}_{} =\bv_1-\bv_2$, $\widetilde{\vec{w}}_{} =\widetilde{\bv}_1-\widetilde{\bv}_2$ and $\overline{\vec{w}}_{} =\overline{\bv}_1-\overline{\bv}_2$
 We subtract the equation for $\bv_1$ from the equation for $\bv_2$ and test it with $\bw$.
 
In the following we distinguish between two cases.\\

\underline{Case 1:} $   2\theta_1+ \theta_2=\frac{1}{2},        0 \le \theta_1 < \frac{1}{4}$ and   $0 < \theta_2< \frac{1}{2}.$ 
 
  We get using successively Cauchy-Schwarz inequality, Young inequality, embedding theorem and the relations (\ref{TKEbis}) and (\ref{TKE}).
 \begin{equation}
 \begin{split}
  \displaystyle 
\|{\vec{w}}_{,t}\|_{2}^{2} +\nu \|\nabla{\vec{w}}_{}\|_{2}^{2} & \le \displaystyle
\frac{4}{\nu} \|\widetilde{\vec{w}}
\overline{\bv}_{1} \|_{2}^{2} +  \frac{4}{\nu} \|\widetilde{\vec{v}_2}
\overline{\vec{w}} \|_{2}^{2}    \\
&\le \displaystyle \frac{4}{\nu} \|\widetilde{\vec{w}}\|_{\frac{1}{2}-\theta_2,2}^{2}
\|\overline{\bv}_{1} \|^{2}_{1+\theta_2}+  \displaystyle\frac{4}{\nu} \|\overline{\vec{w}}\|_{\frac{1}{2}-2\theta_1,2}^{2}
\|\widetilde{\bv}_{2} \|^{2}_{1+2\theta_1}  \\ &\le \displaystyle \frac{1}{\alpha^{{2\theta_1}+{2\theta_2}}}\frac{4}{\nu} \|{\vec{w}}\|_{2}^{2}
\left(\|{\bv}_{1} \|^{2}_{1,2} + 
\|{\bv}_{2} \|^{2}_{1,2}\right).
\end{split}
\end{equation}


\underline{Case 2:} $    \theta_1=0$ and        $ \theta_2 = \frac{1}{2}.$ \\
In this case we have that 
\begin{equation}
\label{vbar11demi}
 \overline{\bv}^N \in L^{\infty}(0,T ; \vec{H}^{1}_{\diver}) \cap L^{2}(0,T ; \vec{H}^{2}_{\diver}),
 \end{equation}
 We get using successively Cauchy-Schwarz inequality, Young inequality, embedding theorem and the relation (\ref{TKE}).

 \begin{equation}
 \begin{split}
  \displaystyle 
\|{\vec{w}}_{,t}\|_{2}^{2} +\nu \|\nabla{\vec{w}}_{}\|_{2}^{2} & \le \displaystyle
\frac{4}{\nu} \|{\vec{w}}
\overline{\bv}_{1} \|_{2}^{2} +  \frac{4}{\nu} \|{\vec{v}_2}
\overline{\vec{w}} \|_{2}^{2}    \\
&\le \displaystyle \frac{4}{\nu} \|{\vec{w}}\|_{2}^{2}
\|\overline{\bv}_{1} \|^{2}_{\frac{3}{2}+\epsilon}+  \displaystyle\frac{4}{\nu} \|\overline{\vec{w}}\|_{1,2}^{2}
\|{\bv}_{2} \|^{2}_{\frac{1}{2},2}  \\ &\le \displaystyle \frac{1}{\alpha^{2}}\frac{4}{\nu} \|{\vec{w}}\|_{2}^{2}
\left(\|{\bv}_{1} \|^{2}_{2,2} + 
\|{\bv}_{2} \|^{2}_{1,2}\right).
\end{split}
\end{equation}
%
%
%
%

Using Gronwall's inequality we conclude 
the continuous dependence of the solutions on the inital data in the $L^{\infty}([0,T],\vec{L}^2_{})$  norm. In particular, if ${\vec{w}}^{}_{0}=0$ then ${\vec{w}}=0$ and the solutions are unique for all $t \in [0,T] .$  Since $T>0$ is arbitrary this solution may be uniquely extended for all time.\\ 
This finish the proof of Theorem \ref{TH1}.\\

\section{Hausdroff dimension of the time singular set in the subcritical case: $   2\theta_1+ \theta_2 < \frac{1}{2} $      }

 The aim in this section is to establish an upper bound for the Hausdorff dimension of the time singular set $\mathcal{S}_{\theta_1,\theta_2}$ of the solutions $\overline{\vec{v}}$ of    \eqref{BM}--\eqref{TKE}, see Theorem \ref{main} below.  We know, thanks  to 
Scheffer's work \cite{S76,S77-2}, that  if $\vec{v}$ is a weak Leray solution of the Navier-Stokes equations then the 
$\frac{1}{2}$-dimensional Hausdorff measure of the time singular set of $\vec{v}$  is zero. Further, 
when $2\theta_1+\theta_2=\frac{1}{2}$, we  proved in the above section  the existence  of a unique  regular weak solution to the model   \eqref{BM}--\eqref{TKE}.
Therefore, it is intersecting to understand how the  time singular set $\mathcal{S}_{\theta_1,\theta_2}(\overline{\vec{v}})$  may depend on the regularization 
parameters $\theta_1$ and $\theta_2$. \\

We divide this section into four subsections.  One is devoted to prove the existence of weak solutions. The second one is devoted to prove the existence of a unique strong solution. An additional subsection is devoted to the defintions of the Hausdorff dimesion and the singular time set. The final subsection is devoted the the proof of Theorem \ref{main} which is the main result of this section.

\subsection{ Existence of weak solutions}

\begin{Theorem}
\label{1deuxieme}
Assume that $2\theta_1 + \theta_2 < \frac{1}{2}$. Let  $\bef \in L^{2}
(0,T;\vec{H}^{-1}_{})$ be a divergence free function and  $\overline{\bv}_0 \in \vec{H}^{\theta_2}_{\diver}$. Then  for any $T>0$ there exist $(\overline{\bv}, \bv,p)$ a  weak distributional solution to \eqref{BM}--\eqref{TKE}  such that
\begin{equation}
\begin{split}
  \overline{\vec{v}_{}} \in  C_{weak}([0,T];\vec{H}^{\theta_2}_{\diver}) \cap L^{2}([0,T];\vec{H}^{1+\theta_2}_{\diver}), \label{ns12}\\
   \vec{v}_{} \in  C_{weak}([0,T];\vec{H}^{-\theta_2}_{\diver}) \cap L^{2}([0,T];\vec{H}^{1-\theta_2}_{\diver}),\\
   \end{split}
   \end{equation}
   \begin{equation}
\begin{split}
\displaystyle \frac{ \partial\overline{\vec{v}}_{}^{}}{ \partial t} \in {L^{\frac{5}{3-2\theta_1-2\theta_2}}([0,T];W^{-1+2\theta_2,\frac{5}{3-2\theta_1-2\theta_2} }(\mathbb{T}_3)^3))},\\
 \displaystyle \frac{ \partial\vec{v}_{}^{}}{ \partial t} \in {L^{\frac{5}{3-2\theta_1-2\theta_2}}([0,T];W^{-1,\frac{5}{3-2\theta_1-2\theta_2} }(\mathbb{T}_3)^3))},
  \end{split}
   \end{equation}
   \begin{equation}
\begin{split}
   {p}_{} \in  L^{\frac{5}{3-2\theta_1-2\theta_2}}([0,T],L^{{\frac{5}{3-2\theta_1-2\theta_2}}}(\mathbb{T}_3)),
   \end{split}
   \end{equation}
\begin{equation}
\begin{split}
\int_0^T \langle \frac{ \partial\vec{v}_{}^{}}{ \partial t} , \bw \rangle  -  (\widetilde{\vec{v}} \otimes \overline{\vec{v}}, \nabla
\bw) +  \nu(
\nabla \bv, \nabla \bw )- (p, \diver \bw) \; dt\\
\qquad 
= \int_0^T \langle \bef, \bw \rangle \; dt
\qquad \textrm{ for all } \bw\in L^{\frac{5}{2+2\theta_1+2\theta_2}}(0,T; W^{1,\frac{5}{2+2\theta_1+2\theta_2}}(\mathbb{T}_3)^3)),
\end{split}\label{weak1001}
\end{equation}
where the velocity $\vec{v}$ verifies 
\begin{equation}\begin{split}
\sup_{t\in(0,T)}\|{\vec{v}}(t)\|_{-1,2}^2+ \nu\int_0^T \|{\vec{v}}(t)\|_{2}^2 dt\le \|{\vec{v}}_0\|_{-1,2}^2 + \int_0^T \langle \bef,\overline{\vec{v}} \rangle \; dt
, \label{regularity}
\end{split}
\end{equation}
or equivalently $\overline{\bv}$ verifies 
\begin{equation}
\begin{array}{lll}
 \label{apriori12lerayinequality}
\displaystyle \sup_{t\in(0,T)}\left(\|\overline{\bv} \|_{2}^2 + \|\overline{\bv}\|_{\theta_2,2}^2\right)+\displaystyle 2\nu\int_{0}^{T}\left(\|\overline{\bv} \|_{1,2}^2  + \|\overline{\bv} \|_{1+\theta_2,2}^2 \right)dt \\
\quad \quad \quad \le \displaystyle 2\int_{0}^{T}\langle \bef, \overline{\bv} \rangle dt + \left(\|\overline{\bv}_0 \|_{2}^2 + \|\overline{\bv}_0\|_{\theta_2,2}^2\right),
\end{array}
\end{equation}
  and the intial data is attained in the following sense 
\begin{equation}
\begin{split}
\label{initial1}
\lim_{t\to 0+}\left(\|\vec{v}(t)-\vec{v}_0\|_{-\theta_2,2}^2 \right)
=0.\\
\lim_{t\to 0+}\left(\|\overline{\vec{v}}(t)-\overline{\vec{v}}_0\|_{\theta_2,2}^2 \right)
=0.
\end{split}
\end{equation}
\end{Theorem}
\textbf{Proof of Theorem \ref{1deuxieme} }
The proof of Theorem \ref{1deuxieme} follows the lines of the proof of the above Theorem the only difference is that $\overline{\bv}$ is not a good test function in the weak formulation and thus by using the weak convergence  we get the inequality (\ref{regularity}) instead of an equality. 

It remains to show weak continuity in \eqref{ns12} and \eqref{initial1}. This is standard for Navier Stokes equation, we refer the reader to  \cite[Lemma~1.4]{Te79} and we omit more  details.

\subsection{Strong solution}
\begin{Theorem}
\label{2deuxieme}
Let $\vec{f} \in L^{2}([0,T],\vec{L}^2_{\diver})$ and $\overline{\vec{v}}_{0} \in \vec{H}^{1-\theta_2}_{\diver}$. Assume that $0\leq 2\theta_1+\theta_2< \frac{1}{2} $.
 Then there exists $T_*:= T_*(\overline{\vec{v}}_{0})$, determined by (\ref{maximal time}), and  there exists a unique strong solution  $ \overline{\vec{v}}$ to  
 \eqref{BM}--\eqref{TKE} on $[0, T_*[$ satisfying:
$$ \overline{\vec{v}_{}} \in  C([0,T_*[;\vec{H}^{\theta_2}) \cap L^{2}([0,T_*[;\vec{H}^{1+\theta_2}_{\diver}),$$  $$ \displaystyle \frac{\partial\vec{v}_{}^{}}{\partial t} \in {L^{2}([0,T_*[;\vec{L}^{2})} \quad \hbox{ and } \quad   
{p}_{} \in  L^{2}([0,T_*[,W^{1,2}(\mathbb{T}_3))\,.$$
\end{Theorem}
\textbf{Proof of Theorem \ref{2deuxieme}} 
Taking the $L^2$-inner product of \eqref{BLM}  with $-\Delta \overline{\vec{v}}$ and integrating by parts. Using  the incompressibility of the velocity field and the duality relation  
combined with H\"{o}lder inequality and Sobolev injection, we obtain
\begin{equation}
\begin{array}{lll}
\displaystyle \frac12\frac{d}{dt}\left(\|\overline{\bv} \|_{1,2}^2 + \|\overline{\bv}\|_{1+\theta_2,2}^2\right)+\displaystyle \nu\left(\|\overline{\bv} \|_{2,2}^2  + \|\overline{\bv} \|_{2+\theta_2,2}^2 \right)\\
 \quad \le  \displaystyle \int_{\tore}  |\widetilde{\vec{v}}\cdot \nabla \overline{\vec{v}} \Delta \overline{\vec{v}}| dx +\displaystyle \int_{\tore}|\vec{f} \cdot \Delta \overline{\vec{v}}dx|\\
\quad \le\displaystyle \alpha^{-2\theta_1-2\theta_2}\|\overline{\vec{v}}\|_{1+\theta_2,2}  \|\overline{\vec{v}}\|_{ \frac{3}{2}-2\theta_1,2}   \|\overline{\vec{v}}\|_{2+ \theta_2,2}+   \|\vec{f}\|_{2} \| \overline{\vec{v}}\|_{2,2}.
\end{array}
\end{equation}
Interpolating between 
$\vec{H}^{1+\theta_2}$ and $\vec{H}^{2+\theta_2}$ we get 
\begin{equation}
\begin{array}{lll}
\displaystyle \frac12\frac{d}{dt}\left(\|\overline{\bv} \|_{1,2}^2 + \|\overline{\bv}\|_{1+\theta_2,2}^2\right)+\displaystyle \nu\left(\|\overline{\bv} \|_{2,2}^2  + \|\overline{\bv} \|_{2+\theta_2,2}^2 \right)\\
\quad \le \displaystyle \alpha^{-2\theta_1-2\theta_2}\|\overline{\vec{v}}\|_{1+\theta_2,2}^{\frac{3}{2}+\theta_2+2\theta_1}   \displaystyle \|\overline{\vec{v}}\|_{2+\theta_2,2}^{\frac{3}{2}- \theta_2 + 2 \theta_1} +   \|\vec{f}\|_{2} \|\overline{\vec{v}}\|_{2,2}.
\end{array}
\end{equation}
Using Young inequality we get
\begin{equation}
\begin{array}{lll}
\label{3.24}
\displaystyle \frac12\frac{d}{dt}\left(\|\overline{\bv} \|_{1,2}^2 + \|\overline{\bv}\|_{1+\theta_2,2}^2\right)+\displaystyle \nu\left(\|\overline{\bv} \|_{2,2}^2  + \|\overline{\bv} \|_{2+\theta_2,2}^2 \right)\\
\quad \le \displaystyle \frac{1}{\nu} \|\vec{f}\|_{2}^2 + \displaystyle C(\alpha, \theta_1, \theta_2 ) \displaystyle\|\nabla \vec{v}\|_{1+\theta_2}^{ \frac{2(3+2\theta_2+4\theta_1)}{1+2\theta_2+ 4\theta_1}}.
\end{array}
\end{equation}

 We get a differential inequality 
 \begin{eqnarray}
\label{differentiel}
 \displaystyle Y^{'} \le \displaystyle C(\alpha, \theta_1, \theta_2, \nu , f ) Y^{\gamma},
 \end{eqnarray}
 where $$ Y(t)= 1+ \|\overline{\vec{v}}\|_{1+\theta_2,2}^2 \hbox{ and }   \gamma = \frac{2(3+2\theta_2+4\theta_1)}{1+2\theta_2+ 4\theta_1}$$
 We conclude that 
 $$ Y(t) \le \frac{Y(0)}{ ( 1-2Y(0)^{\gamma -1} C(\alpha,\theta_1, \theta_2, \nu, \vec{f}  )t  )^{\frac{1}{\gamma-1}}  }$$
 as long as $ \displaystyle t < \frac{1}{2Y(0)^{\gamma -1} C(\alpha, \theta_1, \theta_2, \nu, f )}$,  and  thus we obtain 
  $$ \sup_{t \in [0, T_*[} \|\overline{\vec{v}}\|_{1+\theta_2,2}^2 \le 2(1+ \|\overline{\vec{v}}_{0}\|_{1+\theta_2,2}^2)$$ 
  \begin{eqnarray}
  \label{maximal time}
  \hbox{for } \displaystyle t \le  T_{*}:= \frac{3}{8  C(\alpha, \theta_1, \theta_2, \nu, f )}\frac{1}{(1+\|\overline{\vec{v}_{0}}\|_{1+\theta_2,2}^2)^{\gamma -1} }.
  \end{eqnarray}
Integrating  (\ref{3.24}) with respect to  time on $[0, T_{*}] $  gives the following  estimates
 $$ \int^{T_*}_{0}\| \overline{\vec{v}}(t)\|_{2+\theta_2,2}^2 dt \le M_{}(T_*),$$
 where  $$ M_{}(T_*)= \frac{1}{\nu}\left(  \|\overline{\vec{v}_{0}}\|_{1+\theta_2,2}^2 + \frac{2}{\nu} \int^{T_*}_{0}\|\vec{f}\|_{2}^2 dt +C(\alpha,\theta_1,\theta_2) [ 2(1+ \|\overline{\vec{v}}_{0}\|_{1+\theta_2,2}^2)]^{\gamma}\right). $$

\subsection{The Hausdorff dimension and singular set}
The basic facts about Hausdorff measure can be found in \cite{Federer}. The following defintion can be found 
in \cite{RT83} 
\begin{Def}
\label{definitionhaussdrof}
Let $X$ be a metric and let $a > 0$. The $a$-dimensionnal Hausdorff measure of a subset $Y$ of $X$ is 
$$ 
\displaystyle 
\mu_a(Y)=\lim_{\epsilon \searrow 0}\mu_{a,\epsilon}(Y)=\sup_{\epsilon > 0}\mu_{a,\epsilon}(Y)
$$  
where
$$ 
\displaystyle
\mu_{a,\epsilon}(Y)=\inf \sum_{j}(\diameter B_j)^a,
$$
the infimum being taken over all the coverings of $Y$ by balls $B_j$ such that $\diameter B_j \le \epsilon$.
\end{Def}
\begin{Def}
Let $T>0$. We denote by the time singular set of $\overline{\vec{v}}(t) $, weak solution of \eqref{BM}--\eqref{TKE} given by Theorem \ref{1deuxieme}, the set of $t \in [0, T]$ on wich $\overline{\vec{v}}(t) \notin \vec{H}^{1+\theta_2}(\mathbb{T}_3)^3 $. 
\end{Def}
\subsection{ Dimesion of the time singular set}
The main result of the section is the following theorem.
\begin{Theorem}
\label{main} Let $\overline{\vec{v}}$ be any weak Leray solution to \eqref{BM}--\eqref{TKE} given by  Theorem \ref{1deuxieme} ( We suppose also that the externel force $ \vec{f} \in L^{2}([0,T],\vec{L}^2_{\diver})$ ). 
Then for any $T >0$ the $  \frac{1-2\theta_2 - 4\theta_1}{2} $-dimensional Hausdorff measure of the time singular set of $\vec{v}$   is zero.
\end{Theorem}
\textbf{Proof of Theorem \ref{main}}\\

\textbf{Step 1:}( Structure of the time singularity set)
We begin by the following Lemma that caracterize the structure of the time singularity set of a weak solution to \eqref{BM}--\eqref{TKE}.
\begin{Lem}
\label{temam book}
We assume that $\overline{\vec{v}}_{0} \in \vec{H}^{1+\theta_2}_{\diver}$,  $\vec{f} \in L^{2}([0,T],\vec{L}^2_{\diver})$ and $ \overline{\vec{v}}$ is any weak solution to \eqref{BM}--\eqref{TKE} given by  Theorem \ref{1deuxieme}. Then there exist an open set $ \mathcal{O} $ of $(0,T) $  such that:\\
 (i) For all $t \in   \mathcal{O} $ there exist  $t \in (t_1, t_2)\subseteq (0,T) $ such that 
$ \vec{v} \in C((t_1, t_2), \vec{H}^{1+\theta_2}) $.\\
(ii) The  Lebesgue measure of $[0,T] / \mathcal{O} $  is zero.
\end{Lem}
\textbf{Proof of Lemma \ref{temam book}.} Since $\overline{\vec{v}}\in  C_{weak}([0,T];\vec{H}^{1+\theta_2}) $, $\overline{\vec{v}_{}}(t)$ is well defined for every $t$ and we can define
$$ 
\Sigma=\{ t \in [0,T], \overline{\vec{v}_{}}(t) \in \vec{H}^{1+\theta_2} \},
$$
$$ 
\Sigma^c=\{ t \in [0,T], \overline{\vec{v}_{}}(t) \notin \vec{H}^{1+\theta_2} \},
$$
$$ 
\mathcal{O}=\{ t \in (0,T),  \exists \epsilon > 0, \overline{\vec{v}_{}}  \in  C((t-\epsilon, t+\epsilon), \vec{H}^{1+\theta_2})  \}.
$$
It is clear that $\mathcal{O}$ is open.\\
Since $  \overline{\vec{v}}  \in  L^{2}([0,T];\vec{H}^{1+\theta_2})$, $\Sigma^c$ has Lebesgue measure zero. Let us take $t_0$ such  that  $ t_0  \in \Sigma, $  and $ t_0  \notin \mathcal{O} $, then according to Theorem \ref{2deuxieme}, their exists $ \epsilon >0$ such that $\overline{\vec{v}} \in   C((t_0, t_0+\epsilon), \vec{H}^{1+\theta_2}) $. So, $ t_0$  is the left end of one of the connected components of $\mathcal{O} $. Thus $\Sigma  / \mathcal{O} $  is countable and  $ [0,T] /  \mathcal{O}$  has Lebesgue measure zero. 
This finishes the proof of Lemma \ref{temam book}.\\
\begin{Rem}
We deduce from Theorem \ref{2deuxieme} that, if $(\alpha_i, \beta_i )$ , $i \in I$, is one of  the connected components of $ \mathcal{O}$, then 
$$ \lim_{t\rightarrow \beta_i} \| \overline{\vec{v}}(t)\|_{1+\theta_2,2}=+\infty.$$
Indeed, otherwise Theorem \ref{2deuxieme} would show that there exist an $\epsilon >0$ such that $\overline{\vec{v}} \in  C((\beta_i, \beta_i+\epsilon), \vec{H}^{1+\theta_2})$ and $\beta_i$ would not be the end of an connected  components of $ \mathcal{O}$. 
\end{Rem}

\textbf{Step 2:}(Main estimate)
We have the following Lemma:
\begin{Lem}
\label{temam book2}
Under the same notations of Lemma \ref{temam book}
Let $( \alpha_i, \beta_i) $, $i \in I$, be the connected components of $ \mathcal{O}$. Then
\begin{eqnarray}
  \label{leray reformulation}
\displaystyle  \sum_{i \in I}( \beta_i - \alpha_i)^{\frac{1-2\theta_2 - 4\theta_1}{2}}<\infty
\end{eqnarray}
\end{Lem}
\textbf{Proof of Lemma \ref{temam book2}.} 
Let $( \alpha_i, \beta_i) $ be one of these connected components and let $t \in ( \alpha_i , \beta_i) \subseteq \mathcal{O}$.   Since  
$\overline{\vec{v}} \in  C_{weak}([0,T];\vec{H}^{\theta_2}) \cap L^{2}([0,T];\vec{H}^{1+\theta_2})$, $\overline{\vec{v}}(t)$ is well defined for every $t \in ( \alpha_i , \beta_i)$ and $t$ can be chosen   such that $\overline{\vec{v}}(t) \in \vec{H}^{1+\theta_2} $.  According to Theorem \ref{2deuxieme}, inequality (\ref{maximal time}), and since  $ \| \overline{\vec{v}}(\beta_i)\|_{1+\theta_2}=+\infty$,  for $t \in ( \alpha_i , \beta_i) $ we have, 
  \begin{eqnarray*}
  \label{maximal time2}
  \displaystyle \beta_i - t \ge   \frac{1}{  C(\alpha, \theta_1, \theta_2, \nu, f )}\frac{1}{(1+\|\overline{\vec{v}}(t)\|_{1+\theta_2,2}^2)^{\gamma -1} },
  \end{eqnarray*}
   where we have used that $\gamma= \frac{3+2\theta_2+4\theta_1}{1+2 \theta_2 +4 \theta_1}   >1 $.\\ 
  Thus
  \begin{eqnarray*}
  \label{maximal time3}
  \displaystyle \frac{ C(\alpha,\theta_1, \theta_2, \nu, \vec{f} )}{(\beta_i - t)^{\frac{1}{\gamma -1}} }  \le  1+\|\overline{\vec{v}}(t)\|_{1+\theta_2,2}^2.
  \end{eqnarray*}
  Then we integrate on $( \alpha_i , \beta_i) $ to obtain 
   \begin{eqnarray*}
  \label{maximal time4}
  \displaystyle { C(\alpha, \theta_1, \theta_2,\nu, \vec{f} )}{(\beta_i - \alpha_i)^{\frac{-1}{\gamma -1}+1} }  \le (\beta_i - \alpha_i) +\displaystyle  \int_{\alpha_i}^{\beta_i} \|\overline{\vec{v}_{}}(t)\|_{1+\theta_2,2}^2 dt,
  \end{eqnarray*}
  Adding all these relations for $ i \in I $ we obtain 
   \begin{eqnarray*}
  \label{maximal time5}
  \displaystyle { C(\alpha, \theta_1, \theta_2, \nu, \vec{f} )} \sum_{i \in I}{(\beta_i - \alpha_i)^{\frac{-1}{\gamma -1}+1} }  \le T +\displaystyle  \int_{0}^{T} \|\overline{\vec{v}}(t)\|_{1+\theta_2,2}^2 dt.
  \end{eqnarray*}
 This finishes the proof of Lemma \ref{temam book2}.\\

\textbf{Step 3:}(Recovering argument)
We set $ \mathcal{S}=\mathcal{S}_{\theta_1,\theta_2}(\overline{\vec{v}})=[0,T]\setminus \mathcal{O}$. We have to prove that the $\frac{1-2\theta_2 - 4\theta_1}{2}$-dimensional Hausdorff measure of $\mathcal{S} $ is zero. 
Since the Lebesgue measure of $ \mathcal{O}$ is finite ,i.e.
 \begin{eqnarray}
 \label{finitemeasure}
 \displaystyle \sum_{i \in I}(\beta_i -\alpha_i) < \infty,
  \end{eqnarray}
it follows from Lemma \ref{temam book2} that for every $\epsilon >0$  there exist a finite part $I_{\epsilon} \subset I$ such that
 \begin{eqnarray}
\displaystyle  \sum_{i \in I\setminus I_{\epsilon}}(\beta_i -\alpha_i) \le \epsilon
 \end{eqnarray}
 and
   \begin{eqnarray}
   \label{keyone}
\displaystyle  \sum_{i \in I\setminus I_{\epsilon}}(\beta_i -\alpha_i)^{\frac{1-2\theta_2 - 4\theta_1}{2}} \le \epsilon
 \end{eqnarray}
 Note that $ \mathcal{S} \subset [0,T] \setminus \displaystyle  \bigcup_{i \in I_{\epsilon}}(\alpha_i,\beta_i)$ and the set $  [0,T] \setminus\displaystyle  \bigcup_{i \in I_{\epsilon}}(\alpha_i,\beta_i)$ 
 is the union of finite number of mutually disjoint closed intervals, say $B_j$, for $j=1,...,N$. Our aim now is to  show that the $\diameter B_j \le \epsilon$.
 Since the intervals $(\alpha_i, \beta_i)$ are mutually disjoint, each interval $\displaystyle  (\alpha_i, \beta_i)$, $\displaystyle  i \in I\setminus I_{\epsilon}$, 
 is included in one, and only one, interval $B_j$. We denote by $ I_j$ the set of indice $i$ such that $(\alpha_i, \beta_i ) \subset B_j$. It is clear that 
 $ \displaystyle I_{\epsilon},I_1,...,I_N $ is a partition of $I$ and we have $B_j=(\bigcup_{i \in I_{j}}(\alpha_i,\beta_i))\cup(B_j\cap\mathcal{S})$ for all 
 $j=1,...,N$. It follows from (\ref{finitemeasure}) that
 \begin{eqnarray}
 \label{dimetereepsilon}
 \diameter B_j = \sum_{i \in I_{j}}(\beta_i -\alpha_i) \le \epsilon.
  \end{eqnarray}
 Finally in virtue of  the definition \ref{definitionhaussdrof} and estimates (\ref{dimetereepsilon}), (\ref{keyone})  and since $\displaystyle l^\delta \hookrightarrow l^1$ for all $ 0< \delta < 1$ we have
  \begin{equation}
  \begin{array}{ccccccc}
\displaystyle \mu_{\frac{1-2\theta_2 - 4\theta_1}{2},\epsilon}(\mathcal{S}) &\le& \displaystyle \sum_{j=1}^{N}(\diameter B_j)^{\frac{1-2\theta_2 - 4\theta_1}{2}}\\
\displaystyle &\le& \displaystyle \sum_{j=1}^{N}\left(\sum_{i \in I_{j}}(\beta_i -\alpha_i) \right)^{\frac{1-2\theta_2 - 4\theta_1}{2}}\\
\displaystyle &\le& \displaystyle\sum_{j=1}^{N}\sum_{i \in I_{j}}\left(\beta_i -\alpha_i\right)^{\frac{1-2\theta_2 - 4\theta_1}{2}}\\
\displaystyle &=& \displaystyle \sum_{i \in I\setminus I_{\epsilon}}\left(\beta_i -\alpha_i \right)^{\frac{1-2\theta_2 - 4\theta_1}{2}} \le \epsilon.
\end{array}
 \end{equation}
 Letting $\epsilon \rightarrow 0,$ we find $\displaystyle \mu_{\frac{1-2\theta_2 - 4\theta_1}{2}}(\mathcal{S})=0$ and this  completes the proof.

\end{document}